\documentclass[12pt,pdflatex]{article}
\usepackage[dvipdfmx]{graphicx}
\usepackage{subfig}
\usepackage{pict2e}
\usepackage{cite}
\usepackage[T1]{fontenc}
\usepackage{microtype}
\usepackage{colonequals}
\usepackage{enumerate}

\usepackage{amsmath}
\usepackage{amssymb}
\usepackage{theorem}
\usepackage{color}
\usepackage{latexsym}

\newcommand{\A}{{\mathcal A}}

\newcommand{\Zb}{{\mathbb{Z}}}

\newcommand{\Nmb}{{\mathbf{N}}}
\newcommand{\Nmbp}{{{\mathbf{N}^{\prime}}}}

\newcommand{\qed}{\hfill\hbox{\rule{6pt}{6pt}}}
\newtheorem{theorem}{Theorem}[section]
\newtheorem{lemma}[theorem]{Lemma}

\newtheorem{prop}[theorem]{Proposition}

\newtheorem{definition}[theorem]{Definition}

\theorembodyfont{\rmfamily}\newtheorem{remark}[theorem]{Remark}
\theorembodyfont{\rmfamily}\newtheorem{example}[theorem]{Example}

\DeclareMathOperator{\dist}{dist}

\allowdisplaybreaks

\title{Universal parameterized family of distributions of runs
\footnote{
Pleriminary version has been presented in Probability Symposium 2022 Kyoto Univ.
\cite{takahashi2022ProbSympo}}}
\author{Hayato Takahashi\footnote{Random Data Lab.~Inc., Tokyo 1210062,\ Email: hayato.takahashi@ieee.org}}

\begin{document}
  \maketitle

\begin{abstract}
We present explicit formulae for parameterized families of probabilities of the number of nonoverlapping words 
and  increasing nonoverlapping words
in independent and identically distributed (i.i.d.) finite valued random variables, respectively.
Then we provide an explicit formula for a parameterized family of probabilities of the number of runs, which 
 generalizes \(\mu\)-overlapping probabilities for \(\mu\geq 0\) in i.i.d.~binary valued random variables.
We also demonstrate exact probabilities of the number of runs whose size are exactly given numbers (Mood 1940).
The number of arithmetic operations required to compute our  formula for generalized probabilities of runs
is linear order of sample size for fixed number of parameters and range.
To analyse these number of arithmetic operations for unbounded number of parameters,
we show an asymptotic formula for
the number of integer partitions that are less than or equal to given number as a special case of Meinardus's theorem.
Mathematics Subject Classiﬁcation: 05A15,  62E15, 11P82
\end{abstract}

\section{Introduction}\label{secintro}
We study probabilities  of the number of words in finite valued i.i.d.~random variables (probabilities of words for short).
The probabilities of words play important role in statistics, DNA analysis, information theory
 \cite{{koutrasBook},{bertheRigo},{feller1},{JaquetandSzpankowski},{lothire},{Mood40},{bookRobin},
{Wald40},{waterman},{wolf88}}.

Generating functions of the probabilities of words obtained by inductive relations of words on sample size are rational functions
\cite{{Bassinoetal2011}, {blomThorburn}, {chrysaphinou88}, {Flajolet2009},{Goulden83}, 
 {Guibas81},{RegnierSpankowski98}}.
 However, unless all poles of rational function are known, we do not have its partial fraction expansion cf.~Chapter 11 Section 4 pp.~275 \cite{feller1}.
In \cite{{feller1},{JaquetandSzpankowski},{bookRobin}} 
approximations and recurrence formulae for the probabilities of words are given.
In \cite{{uppuuri},{antzou99}} explicit formulae are obtained by directly expanding rational generating functions into power series, i.e.,
\(\frac{f}{1-g}=f\sum_n g^n\) if \(|g|<1\) for polynomials \(f\) and \(g\).

A word that consists of the same letter is called a {\it run}. 
The number of runs depends on the counting manner. 
Let \(0^m\) be the word that consists of \(m\) zeros. 
For \(x\in\{0,1\}^n\), let\\
(i) \(E_{n,m}(x)\), the number of \(0^m\) of size exactly \(m\) in \(x\) \cite{{Mood40},{Fu94}},\\
(ii) \(G_{n,m}(x)\), the number of \(0^m\) of size greater than or equal to \(m\) in \(x\) \cite{{Fu94},{antzou99}},\\
(iii) \(N_{n,m}(x)\), the number of nonoverlapping  \(0^m\) in \(x\) \cite{{Godbole90},{Hirano86},{Muselli96},{antzou99},{Fu94},{feller1}},\\
(iv) \(M_{n,m}(x)\), the number of overlapping  \(0^m\) in \(x\) \cite{{ling88},{antzou99},{koutras97},{Fu94},{Godbole92}}, \\
(v) \(L_n(x)\), the size of the longest  run of 0s in \(x\) \cite{{Makri2007},{Phillipou86},{antzou99},{Fu94}}, \\
(vi) \(T_k(x)\), the stopping time \(t\)  such that \(0^k\) first appear in \(x=x_1\cdots x_t\) \cite{{hirano84},{phili83},{uppuuri}}, and\\
(vii) \(N_{n,m,\mu}(x)\),  the \(\mu\)-letters overlapping enumeration in the string \(x\)
\cite{{akihirano2000},{hanAki},{Makri2015}}.

In this paper, we present explicit formulae for parametric exact probabilities of these statistics.
To avoid the difficulty of enumerating overlapping words in Theorem~\ref{msjOkayama}, 
we study probabilities of increasing nonoverlapping words and their finite dimensional generating functions. 
Combining Theorem~\ref{msjOkayama} with Lemma~\ref{mainlemma}, in Theorem~\ref{mainth2}, we
 derive explicit formulae for parameterized probabilities of runs including those of the statistics (i)--(vii) above by a unified manner 
  for binary valued i.i.d.~random variables.

The rest of the paper consists as follows.
In Section~\ref{Secnonoverlapwords} Theorem~\ref{MAINTHA} and \ref{th-moments}, we show explicit formulae for 
parameterized families of joint probabilities of nonoverlapping words 
and their moments for finite valued i.i.d.~radnom variables.
In Section~\ref{seccomplexity}, we study algorithm and complexity to compute our formulae.
We demonstrate the number of arithmetic operations required to compute our  formula for generalized probabilities of runs
is linear order of sample size for fixed number of parameters and range. 
To analyse these number of arithmetic operations for unbounded number of parameters, in Lemma~\ref{eqnMeinardus},
we show an asymptotic formula for
the number of integer partitions that are less than or equal to given number as a special case of Meinardus's theorem (Chap.~6 \cite{andrews}).
In Section~\ref{Secdist}, we study distance among the distributions of runs. 
\section{Joint probabilities of nonoverlapping words}\label{Secnonoverlapwords}
A finite string of a finite alphabet \(\A\) is called a word.
Let \(|x|\) be the length of a word \(x\). 
The word \(xy\) is the concatenation of two words \(x\) and \(y\).
The word \(x^k\) is the \(k\)-times concatenations of a word \(x\), e.g.~\(x^2=xx\).
A word \(x\) is called overlapping if there is a word \(z\) such that \(x\) appears at least 2 times in \(z\) and \(|z|<2|x|\); otherwise
\(x\) is called nonoverlapping. 
A pair of words \((x,y)\in S^2, x\ne y\) is called pairwise overlapping if there is a word \(z\) such that \(x\) and \(y\) appear in \(z\) and \(|z| < |x|+|y|\);
otherwise the pair is called pairwise nonoverlapping.
A finite set of words \(S\) is called nonoverlapping if every \(x\in S\) is nonoverlapping and every pair \((x,y)\in S^2, x\ne y\) is pairwise nonoverlapping;
otherwise \(S\) is called overlapping.
For example, \(\{10\}, \{1100, 10100\},\text{ and } \{00111,00101\}\) are nonoverlapping; \(\{00\}\), \(\{10,01\}\), and \(\{00,11\}\) are overlapping.

In the following, let 
 \(\Nmb(w_1,\ldots,w_h; x_1^n)\) be the number of  words \(w_1,\ldots,w_h\) in an arbitrary position of \(x_1^n\in \A^n\),
 i.e. 
 \[\Nmb(w_1,\ldots,w_h; x_1^n):= 
 (\sum_{i=1}^{n-|w_1|+1} I_{w_1}(x_i^n),\ldots,\sum_{i=1}^{n-|w_h|+1} I_{w_h}(x_i^n)),\] where
 \(x_i^n=x_i\cdots x_n\) and 
 \(I_{w_j}(x_i^n)=1\) if \(x_i\cdots x_{i+|w_j|-1}=w_j\) else 0   for all \(i,j\).
For \(a_1+\cdots+a_l\leq n\), let
\[\dbinom{n}{a_1,\ldots,a_l}=\frac{n!}{a_1!\cdots a_l! (n-\sum a_i)!},
\]
where \(0!=1\).
Let \(P\) be a probability on \(\A\), i.e., \(0\leq P(a)\leq 1\) for \(a\in\A\) and \(\sum_{a\in\A} P(a)=1\).
Set \(P(w)= \prod P(a_i)\) for \(w=a_1\cdots a_{|w|},\ a_i\in\A\).
For example \(P(w)=2^{-|w|}\) for all \(w\) if \(P(0)=P(1)=1/2\) for \(\A=\{0,1\}\).
\begin{theorem}\label{MAINTHA}
Let \(\A\) be a finite alphabet and \(P\)  a probability on \(\A\).
Let \(X_1^n:=X_1X_2\cdots X_n\)  be \(\A\)-valued i.i.d.~random variables from \(P(X_i=a)=P(a)\) for \(a\in\A\).
Let \(w_1,\ldots, w_h\) be nonoverlapping.
Let 
\begin{align}
&A(k_1,\ldots,k_h)=\dbinom{n-\sum_i |w_i| k_i + \sum_i k_i}{k_1,\ldots, k_h}\prod_{i=1}^h P^{k_i}(w_i),\nonumber\\
&B(k_1,\ldots,k_h)=P(\Nmb(w_1,\ldots,w_h; X_1^n)=(k_1,\ldots,k_h)),\nonumber\\
&F_A(z_1,\ldots,z_h)=\sum_{k_1,\ldots,k_h}A(k_1,\ldots,k_h)z^{k_1}\cdots z^{k_h},\text{ and}\nonumber\\
&F_B(z_1,\ldots,z_h)=\sum_{k_1,\ldots,k_h}B(k_1,\ldots,k_h)z^{k_1}\cdots z^{k_h}.\nonumber
\end{align}
 Then
\begin{align}
& A(k_1,\ldots,k_h)=\sum B(t_1,\ldots,t_h)\dbinom{t_1}{k_1}\cdots\dbinom{t_h}{k_h},\nonumber\\
&F_A(z_1,z_2,\ldots,z_h)=F_B(z_1+1, z_2+1,\ldots,z_h+1),\text{ and}\nonumber\\
&P(N(w_1,\ldots,w_h;X_1^n)=(s_1,\ldots,s_h))\nonumber\\
& =\sum_{\substack{k_1,\ldots,k_h\colon\\ \ s_1\leq k_1,\ldots,s_h\leq k_h\\ \sum_i |w_i| k_i\leq n}} (-1)^{\sum_i k_i-s_i}
\dbinom{n-\sum_i |w_i| k_i + \sum_i k_i}{s_1,\ldots, s_h,k_1-s_1,\ldots k_h-s_h}\prod_{i=1}^h P^{k_i}(w_i).\label{eqmain}
\end{align}
\end{theorem}
Proof)
We prove the theorem for \(h=2\). The proof for the general case is similar.
Let \(x\) be a letter such that \(x\notin \A\).
The number of possible allocations such that \(w_1\) and \(w_2\) appear \(k_1\) and \(k_2\) times respectively without overlapping in \(x^n\)  is
\[\dbinom{n-|w_1|k_1-|w_2|k_2+k_1+k_2}{k_1,k_2}.\]
This is because if we replace \(w_1\) and \(w_2\) with additional extra  symbols \(\alpha\) and \(\beta\) in \(x^n\) then the problem  reduces to 
choosing  \(k_1\)  \(\alpha\)s  and  \(k_2\)  \(\beta\)s among the strings of length  \(n-|w_1|k_1-|w_2|k_2+k_1+k_2\).
Let
\begin{equation}\label{basicFuncA}
A(k_1,k_2)\colonequals \dbinom{n-|w_1|k_1-|w_2|k_2+k_1+k_2}{k_1,k_2}P^{k_1}(w_1)P^{k_2}(w_2).
\end{equation}
The function \(A\) is not the probability of   \(k_1\) \(w_1\)s and \(k_2\) \(w_2\)s occurrences in strings of length \(n\), since if
we allow any letters in the remaining place except for \(k_1\) \(w_1\)s and \(k_2\) \(w_2\)s in the string, that string may include extra \(w_1\)s and \(w_2\)s.
Let \(B(t_1,t_2)\) be the probability that \(w_1\) and \(w_2\) appear \(k_1\) and \(k_2\) times, respectively.
We have the following identity,
\[
A(k_1,k_2) =\sum_{k_1\leq t_1, k_2\leq t_2} B(t_1,t_2) \dbinom{t_1}{k_1} \dbinom{t_2}{k_2}.
\]
Then
\begin{align*}
F_A(z_1,z_2) & =\sum_{k_1,k_2} z_1^{k_1} z_2^{k_2} \sum_{k_1\leq t_1, k_2\leq t_2} B(t_1,t_2) \dbinom{t_1}{k_1}\dbinom{t_2}{k_2}  \\
& = \sum_{t_1,t_2} B(t_1,t_2) \sum_{k_1\leq t_1, k_2\leq t_2} \dbinom{t_1}{k_1}\dbinom{t_2}{k_2} z_1^{k_1} z_2^{k_2} \\
&= \sum_{t_1,t_2} B(t_1,t_2) (z_1+1)^{t_1} (z_2+1)^{t_2}\\
& = F_B(z_1+1,z_2+1).
\end{align*}
We have
\begin{align*}
&F_B(z_1,z_2) =F_A(z_1-1,z_2-1)\\
&=\sum_{k_1,k_2\colon |w_1|k_1+|w_2|k_2\leq n} \dbinom{n-|w_1|k_1-|w_2|k_2+k_1+k_2}{k_1,k_2} (z_1-1)^{k_1} (z_2-1)^{k_2} P^{k_1}(w_1) P^{k_2}(w_2)\\ 
&=\sum_{\substack{k_1,k_2,t_1,t_2\colon |w_1|k_1+|w_2|k_2\leq n\\\phantom{k_1,k_2,t_1,t_2}t_1\leq k_1, t_2\leq k_2}}  
\dbinom{n-|w_1|k_1-|w_2|k_2+k_1+k_2}{k_1,k_2}\dbinom{k_1}{t_1}\dbinom{k_2}{t_2}z_1^{t_1}z_2^{t_2}(-1)^{k_1+k_2-t_1-t_2} \\
&\phantom{\sum_{\substack{k_1,k_2,t_1,t_2\colon |w_1|k_1+|w_2|k_2\leq n\\\phantom{k_1,k_2,t_1,t_2}t_1\leq k_1, t_2\leq k_2}}  }
\times P^{k_1}(w_1) P^{k_2}(w_2)\\
&=\sum_{t_1,t_2} z_1^{t_1}z_2^{t_2}\sum_{\substack{k_1,k_2\colon |w_1|k_1+|w_2|k_2\leq n\\ t_1\leq k_1, t_2\leq k_2}}  
(-1)^{k_1+k_2-t_1-t_2}  \dbinom{n-|w_1|k_1-|w_2|k_2+k_1+k_2}{t_1,t_2,k_1-t_1,k_2-t_2} P^{k_1}(w_1) P^{k_2}(w_2).
\end{align*}
and  (\ref{eqmain}) .
\qed

 In \cite{RegnierSpankowski98},  expectation, variance, and central limit theorems for the occurrences of words are shown.
In \cite{rukhin08}, chi-squared tests with nonoverlapping words are studied. 
We give all orders of moments for nonoverlapping words.
Let \(A_{t,s}\colonequals\sum_r \dbinom{s}{r}r^t(-1)^{s-r}\) for all \(t,s=1,2,\ldots\).
Then \(A_{t,s}\) is the number of surjective functions from \(\{1,2,\ldots,t\}\to\{1,2,\ldots,s\}\) for all \(t,s\), see pp.100 Problem 1  \cite{Riordan}. Let \(\lfloor x \rfloor\) be the greatest integer less than or equal to \(x\).
\begin{theorem}\label{th-moments}
Let \(w\) be nonoverlapping. Under the same assumption with \(h=1\) in  Theorem~\ref{MAINTHA},
\[E(N^t(w;X^n))=\sum_{s=1}^{\min \{\lfloor \frac{n}{|w|}\rfloor,t\}} A_{t,s}\dbinom{n-s|w|+s}{s}P^s(w)\]
for all \(t=1,2,\ldots\).
\end{theorem}
Proof)
Let
\(Y_i:=I_{X_i^{i+|w|-1}=w}\) for \(i+|w|-1\leq n\) and \(S(Y_i):=\{i,i+1,\ldots,i+|w|-1\}\).
Since \(w\) is nonoverlapping, we have 
\begin{equation}\label{expect}
E(Y_iY_j)=
\begin{cases}
P(w) &\text{ if }i=j,\\
P^2(w)& \text{ if }S(Y_i)\text{ and }S(Y_j)\text{ are disjoint,}\\
0 &\text{ else.}
\end{cases}
\end{equation}

Let \(Y_{j,i}=Y_i\) for all \(1\leq j\leq t\).
Then  
\begin{align}
E(N^t(w;X^n))&=E((\sum_i Y_i)^t)=E(\prod_{j=1}^t\sum_i Y_{j,i})=E(\sum_{n(1),\ldots,n(t)}\prod_{j=1}^t Y_{j,n(j)}).\label{subeq}
\end{align}
By (\ref{expect}), \(E( \prod_{j=1}^t Y_{j,n(j)})=P^s(w)\) if and only if 
\(\{S(Y_{1,n(1)}),\ldots,S(Y_{t,n(t)})\}=\{S(Y_{l(1)}),\ldots,S(Y_{l(s)})\}\) and \(S(Y_{l(1)}),\ldots,S(Y_{l(s)})\) are disjoint.

The number of possible combination of \(s\) disjoint \(\{S(Y_{l(1)}),\ldots,S(Y_{l(s)})\}\) is \(\dbinom{n-s|w|+s}{s}\).
If \(n<s|w|\) then there is no \(s\) disjoint \(S(Y_i)\)s.
For each disjoint \(\{S(Y_{l(1)}),\ldots,S(Y_{l(s)})\}\),
the number of possible combination of \(n_1,\ldots, n_t\) such that 
\(\{S(Y_{1,n(1)}),\ldots,S(Y_{t,n(t)})\}=\{S(Y_{l(1)}),\ldots,S(Y_{l(s)})\}\) is
\(A_{t,s}\). 
By (\ref{subeq}), we have the theorem.
\qed

\section{Explicit formulae for distributions of runs}\label{Secruns}
First we show  probability functions for increasing nonoverlapping words. 
 Let
\[
\Nmb^{\prime}(w_1,\ldots,w_h; x)\colonequals  (s_1-s_2,s_2-s_3,\ldots,s_h)\text{ if }\Nmb(w_1,\ldots,w_h;x)=(s_1,\ldots,s_h).
\]
For example \(\Nmb(100,1000;1010001)=(1,1)\) and \(\Nmb^\prime(100,1000;1010001)=(0,1)\).
We write \(x\sqsubset y\) if \(x\) is a prefix of \(y\) and \(x\ne y\). For example \(10 \sqsubset 100\).
If  \(w_1\sqsubset w_2\)
 and \((k_1,k_2)=\Nmb(w_1,w_2;x)\) then \(k_1\geq k_2\) for all \(x\).
\begin{definition}
Let
\begin{equation}\label{defCnfun}
C_{n,(w_1,\ldots,w_h)}(x):=t\text{ if }\sum ik_i=t\text{ and }\Nmb^{\prime}(w_1,\ldots,w_h; x_1^n)=(k_1,k_2,\ldots,k_h)\text{ for }
|x|=n,
\end{equation}
where \(w_1\sqsubset w_2\cdots \sqsubset w_h\) be  increasing nonoverlapping words.
\end{definition}
\begin{theorem}\label{msjOkayama}
Let \(\A\) be a finite alphabet and \(P\)  a probability on \(\A\).
Let \(X_1,X_2,\ldots,\) be \(\A\)-valued i.i.d.~finite valued random variables
from \(P(X_i=a)=P(a)\) for \(a\in\A\).
Let \(w_1\sqsubset w_2\cdots \sqsubset w_h\) be  increasing nonoverlapping words and 
\begin{align}
&A(k_1,\ldots,k_h):= \dbinom{n-\sum_i |w_i| k_i + \sum_i k_i}{k_1,\ldots,k_h}\prod_{i=1}^h P^{k_i}(w_i),\nonumber\\
&B(k_1,\ldots,k_h):=  P(\Nmb^{\prime}(w_1,\ldots,w_h; X_1^n)=(k_1,k_2,\ldots,k_h)),\nonumber\\
&F_A(z_1,\ldots,k_h):=  \sum_{\substack{k_1,\ldots,k_h: \\ \sum_i |w_i| k_i\leq n}}A(k_1,\ldots,k_h)z^{k_1}\cdots z^{k_h},\text{ and}\nonumber\\
&F_B(z_1,\ldots,z_h):=  \sum_{\substack{k_1,\ldots,k_h: \\ \sum_i |w_i| k_i\leq n}}B(k_1,\ldots,k_h)z^{k_1}\cdots z^{k_h}.
\end{align}
Then
\begin{align}
&F_A(z_1,\ldots,z_h)=F_B(z_1+1,z_1+z_2+1,\ldots,\sum_i z_i +1)\text{ and}
\label{msjOkayamaeq}\\
&P( C_{n,(w_1,\ldots,w_h)}(X_1^n)=t)
=\sum_{\substack{k_1,\ldots,k_h:\\ \sum |w_i|k_i\leq n,\\ \sum (i-1)k_i\leq t\leq\sum ik_i}}
(-1)^{\sum ik_i-t}\dbinom{n-\sum |w_i|k_i+\sum k_i }{k_1,\ldots,k_h} \dbinom{\sum k_i}{\sum ik_i-t}\prod P^{k_i}(w_i).\label{subexplicitform}
\end{align}
\end{theorem}
Proof)
We show (\ref{msjOkayamaeq}) for \(h=2\). The proof of the general case is similar.
Observe that 
\begin{align}\label{eq-proof-A}
A(k_1,k_2)=&\sum_{k_2\leq t_2, \ k_1+k_2\leq t_1+t_2}  B(t_1,t_2)\dbinom{t_2}{k_2} 
\sum_{0\leq s\leq t_2-k_2}\dbinom{t_2-k_2}{s}\dbinom{t_1}{k_1-s}.
\end{align}
Then
\begin{align*}
F_A(z_1,z_2)=&\sum_{k_1,k_2}z_1^{k_1}z_2^{k_2} \sum_{k_2\leq t_2, \ k_1+k_2\leq t_1+t_2} B(t_1,t_2)
\dbinom{t_2}{k_2}\sum_{0\leq s\leq t_2-k_2}\dbinom{t_2-k_2}{s}\dbinom{t_1}{k_1-s}\\
=&\sum_{t_1,t_2}B(t_1,t_2)\sum_{k_2\leq t_2}\dbinom{t_2}{k_2}z_2^{k_2}  
\sum_{0\leq s\leq t_2-k_2,\ 0\leq k_1-s\leq t_1}\dbinom{t_2-k_2}{s}\dbinom{t_1}{k_1-s}z_1^{k_1}\\
=&\sum_{t_1,t_2}B(t_1,t_2)\sum_{k_2\leq t_2} \dbinom{t_2}{k_2}z_2^{k_2} (z_1+1)^{t_1+t_2-k_2}\\
=&\sum_{t_1,t_2}B(t_1,t_2)(z_1+1)^{t_1+t_2}(\frac{z_2}{z_1+1}+1)^{t_2}\\
=&F_B(z_1+1,z_1+z_2+1).
\end{align*}
Next, set \(z_1=X, z_2=X(X+1),\ldots, z_h=X(X+1)^{h-1}\) in (\ref{msjOkayamaeq}).
Then
\begin{equation}\label{msjOkayamaB}
F_A(X,X(X+1),\ldots,X(X+1)^{h-1})=F_B(X+1,(X+1)^2,\ldots,(X+1)^h).
\end{equation}
By setting \(Y=X+1\) in (\ref{msjOkayamaB}), we have 
\begin{equation}\label{msjOkayamaRev}
F_A(Y-1,(Y-1)Y,\ldots,(Y-1)Y^{h-1})=F_B(Y,Y^2,\ldots,Y^h).
\end{equation}
Since
\[F_B(Y,Y^2,\ldots,Y^h)=\sum_{\substack{k_1,\ldots,k_h: \\ \sum_i |w_i| k_i\leq n}}B(k_1,\ldots,k_h)Y^{\sum ik_i},\]
 \(P(\sum ik_i=t)\) is the coefficient of \(Y^t\) in \(F_B\).
On the other hand, by expanding the left-hand-side of (\ref{msjOkayamaRev}), we have
\begin{align*}
&F_A(Y-1,(Y-1)Y,\ldots,(Y-1)Y^{h-1})\nonumber\\
&=\sum_{k_1,\ldots,k_h} \dbinom{n-\sum |w_i|k_i+\sum k_i}{k_1,\ldots,k_h}
 (Y-1)^{\sum k_i}\prod Y^{(i-1)k_i}P^{k_i}(w_i)\nonumber\\
&=\sum_{k_1,\ldots,k_h} \dbinom{n-\sum |w_i|k_i+\sum k_i}{k_1,\ldots,k_h}\prod P^{k_i}(w_i)
\sum_r \dbinom{\sum k_i}{r} (-1)^rY^{\sum ik_i-r}\nonumber.
\end{align*}
Let \(t=\sum ik_i-r\) then \(0\leq r\leq \sum k_i\Leftrightarrow \sum (i-1)k_k\leq t\leq \sum ik_i\), and we have (\ref{subexplicitform}).
\qed\\

To derive a universal formula for probability functions of runs, we introduce a statistics that represents various types of runs.
\begin{definition}
For \(x\in\{0,1\}^n\), let
\begin{equation}
D_{n,(m_1,\ldots,m_h)}(x)\colonequals t\text{ if }\sum ik_i=t\text{ and }\Nmb^{\prime}(10^{m_1}
\ldots,10^{m_h}; 1x)=(k_1,k_2,\ldots,k_h),
\end{equation}
where \(m_1<\ldots < m_h\).
\end{definition}

\begin{example}
Consider a run \(0^3\) and let \(x=0\ 0\ 0\ 1\ 0\ 0\ 0\ 0\ 0\ 0\ 0\ 1\ 0\ 0\ 0\ 0\).\\
1.~Let \(m_i=0^{3i}\) for \(1\leq i\leq 5\).
Then \(\Nmbp(10^3,10^6,\ldots,10^{15};1x)=(2,1,0,\ldots,0)\) and
\(D_{16,(3,6,\ldots,15)}(x)=\sum ik_i=2+2\cdot 1=4=N_{16,3}(x)\) (0-overlapping enumeration).\\
2.~Let \(m_i=0^{3+2(i-1)}=0^{2i+1}\) for \(1\leq i\leq 7\).
Then \(\Nmbp(10^3,10^5,\ldots,10^{15};1x)=(2,0,1,0,\ldots,0)\) and 
\(D_{16,(3,5,\ldots,15)}(x)=\sum ik_i=2+3\cdot 1=5\) (1-overlapping enumeration).\\
3.~Let \(m_i=0^{3+i-1}=0^{2+i}\) for \(i=1,2,\ldots, 14\) and 
Then \(\Nmbp(10^3,10^4,\ldots,10^{16};1x)=(1,2,0,0,1,0,\ldots,0)\) and
\(D_{16,(3,4,\ldots,16)}(x)=\sum ik_i=1+2\cdot 2+5\cdot 1=8=M_{16,3}(x)\) (2-overlapping enumeration).\\
4.~Let \(m_1=0^3\). Then  \(\Nmbp(10^3;1x)=(3)\) and 
\(D_{16,(3)}(x)=3=G_{16,3}(x)\).
\end{example}
When \(w_i=10^{m_i}\), the difference between \(D_n\) and \(C_n\) is that 
\(D_n\) count \(0^m\) for \(m\geq m_1\) from the beginning of \(x\) while 
\(C_n\) does not.
\begin{lemma}\label{mainlemma}
Let \(X_1,X_2,\ldots,\) be i.i.d.~binary random variables
from \(P(X_i=1)=q\) and \(P(X_i=0)=p\) for all \(i\). 
Let \(m_1<\ldots < m_h\) and \(w_i=10^{m_i}\) for \(1\leq i\leq h\).
Then for all \(t\geq 0\),
\[P(D_{n,(m_1,\ldots,m_h)}(X_1^n)=t)=(P(C_{n+1,(w_1,\ldots,w_h)}(X_1^{n+1})=t)-pP(C_{n,(w_1,\ldots,w_h)}(X_1^n)=t))/q.\]
\end{lemma}
Proof)
Observe that 
\begin{align*}
&\{0x\mid C_{n+1,(w_1,\ldots,w_h)}(0x)=t, |x|=n\}=\{0x\mid C_{n,(w_1,\ldots,w_h)}(x)=t, |x|=n\},\text{ and}\\
&\{1x\mid C_{n+1,(w_1,\ldots,w_h)}(1x)=t, |x|=n\}=\{1x\mid D_{n,(m_1,\ldots,m_h)}(x)=t, |x|=n\}.
\end{align*}
We have
\begin{align*}
P(C_{n+1,(w_1,\dots,w_h)}(X_1^{n+1})=t)&=P(C_{n+1,(w_1,\dots,w_h)}(0X_1^{n})=t)+ P(C_{n+1,(w_1,\dots,w_h)}(1X_1^{n})=t)\\
&=p P(C_{n,(w_1,\dots,w_h)}(X_1^{n})=t)+qP(D_{n,(m_1,\dots,m_h)}(X_1^{n})=t).
\end{align*}
\qed
\begin{theorem}[main theorem]\label{mainth2}
Let \(X_1,X_2,\ldots,\) be i.i.d.~binary random variables
from \(P(X_i=1)=q\) and \(P(X_i=0)=p\) for all \(i\). 
Let \(m_1<\ldots < m_h\) and \(w_i=10^{m_i}\) for \(1\leq i\leq h\).
Then for all \(t\geq 0\),
\begin{align*}
&1.~P( D_{n,(m_1,\ldots,m_h)}(X_1^n)=t)
=\sum_{\substack{k_1,\ldots,k_h:\\ \sum (m_i+1)k_i\leq n+1,\\ \sum (i-1)k_i\leq t\leq \sum ik_i}}
(-1)^{\sum ik_i-t}\dbinom{n+1-\sum m_ik_i }{k_1,\ldots,k_h} \dbinom{\sum k_i}{\sum ik_i-t} q^{k_i-1}p^{k_i\sum m_i}\\
&\phantom{P( D_{n,(m_1,\ldots,m_h)}(X_1^n)=t)=}
-\sum_{\substack{k_1,\ldots,k_h:\\ \sum (m_i+1)k_i\leq n,\\ \sum (i-1)k_i\leq t\leq \sum ik_i}}
(-1)^{\sum ik_i-t}\dbinom{n-\sum m_ik_i}{k_1,\ldots,k_h} \dbinom{\sum k_i}{\sum ik_i-t} q^{k_i-1}p^{1+k_i\sum m_i},\\
&2.~P(N_{n,m,\mu}(X_1^n)=t)=P( D_{n,(m_1,\ldots,m_h)}(X_1^n)=t)\text{ for all }0\leq\mu\leq m-1,\text{ where}\\
&h=\lfloor \frac{n-\mu}{m-\mu}\rfloor\text{ and }m_i=mi-\mu (i-1)\text{ for }1\leq i\leq h,\\
&3.~P(G_{n,m}(X_1^n)=t)=P(D_{n,(m)}=t)=\\
&\sum_{ t\leq k\leq \lfloor\frac{n+1}{m+1}\rfloor} (-1)^{k-t}\binom{n+1-mk}{t,k-t}q^{k-1}p^{km}
-\sum_{ t\leq k\leq \lfloor\frac{n}{m+1}\rfloor} (-1)^{k-t}\binom{n-mk}{t,k-t}q^{k-1}p^{km+1},\\
&4.~P(T_m>n)=P(L_n<m)=P(D_{n,(m)}=0)=\\
&\sum_{0\leq k\leq \lfloor \frac{n+1}{m+1}\rfloor}(-1)^k\dbinom{n+1-mk}{k}q^{k-1}p^{km}
-\sum_{0\leq k\leq \lfloor \frac{n}{m+1}\rfloor}(-1)^k\dbinom{n-mk}{k}q^{k-1}p^{km+1},\text{ and}\\
&5.~P(E_{n,m}(X_1^n)=t)=\\
&\sum_{\substack{k_1,k_2:\\(m+1)k_1+(m+2)k_2\leq n+1,\\ t\leq k_1+k_2}} 
(-1)^{k_1-t}\dbinom{n+1-mk_1-(m+1)k_2}{k_1,k_2}\dbinom{k_1+k_2}{t} 
q^{k_1+k_2-1}p^{k_1m+k_2(m+1)}\nonumber\\
&-\sum_{\substack{k_1,k_2:\\(m+1)k_1+(m+2)k_2\leq n,\\ t\leq k_1+k_2}} 
(-1)^{k_1-t}\dbinom{n-mk_1-(m+1)k_2}{k_1,k_2}\dbinom{k_1+k_2}{t} q^{k_1+k_2-1}p^{k_1m+k_2(m+1)+1}.
\end{align*}
\end{theorem}
Proof)
Part 1 follows from Theorem~\ref{msjOkayama} and Lemma~\ref{mainlemma}.
Part 2 follows from part 1.
Part 3 follows from \(P(G_{n,m}=t)=P(D_{n,(m)}=t)\).
Part 4 follows from \(P(T_m>n)=P(L_n<m)=P(G_{n,m}=0)\).\\
Proof of part 5.
Let \(h=2\), \(w_1=10^m\), and \(w_2=10^{m+1}\) in Theorem~\ref{msjOkayama}. By (\ref{msjOkayamaeq}), we have
\begin{equation}\label{subEq1}
F_A(z_1,z_2)=F_B(z_1+1,z_1+z_2+1).
\end{equation}
Set \(z_1=x-1\) and  \(z_2=1-x\).  We have
\begin{align}
F_A(x-1,1-x)&=F_B(x,1)\nonumber\\
&=\sum_{k_1,k_2: (m+1)k_1+(m+2)k_2\leq n} P(\Nmb^\prime(w_1,w_2)=(k_1,k_2))x^{k_1}\nonumber\\
&=\sum_{k_1}\sum_{k_2: (m+1)k_1+(m+2)k_2\leq n}P(\Nmb^\prime(w_1,w_2)=(k_1,k_2))x^{k_1}\nonumber\\
&=\sum_{k_1} P(\bar{E}_{n,m}=k_1)x^{k_1},\label{subeq1A}
\end{align}
where 
\(\bar{E}_{n,m}(x)\), the number of \(10^m\) of size exactly \(m+1\) in \(x\).

On the other hand,
\begin{align}
F_A(x-1,1-x)&
=\sum_{\substack{k_1,k_2:\\(m+1)k_1+(m+2)k_2\leq n}}
\dbinom{n-(m+1)k_1-(m+2)k_2+k_1+k_2}{k_1,k_2}P^{k_1}(w_1)P^{k_2}(w_2)\nonumber\\
&\phantom{\sum_{\substack{k_1,k_2:\\(m+1)k_1+(m+2)k_2\leq n}}}\times (x-1)^{k_1}(1-x)^{k_2}\nonumber\\
&=\sum_{\substack{k_1,k_2:\\(m+1)k_1+(m+2)k_2\leq n}}(-1)^{k_2}\dbinom{n-(m+1)k_1-(m+2)k_2+k_1+k_2}{k_1,k_2}P^{k_1}(w_1)P^{k_2}(w_2)\nonumber\\
&\phantom{\sum_{\substack{k_1,k_2:\\(m+1)k_1+(m+2)k_2\leq n}}}\times (x-1)^{k_1+k_2}\nonumber\\
&=\sum_{\substack{k_1,k_2,t:\\(m+1)k_1+(m+2)k_2\leq n\\ t\leq k_1+k_2}}
(-1)^{k_1+2k_2-t}\dbinom{n-(m+1)k_1-(m+2)k_2+k_1+k_2}{k_1,k_2}\dbinom{k_1+k_2}{t}\nonumber\\
&\phantom{\sum_{\substack{k_1,k_2,r:\\(m+1)k_1+(m+2)k_2\leq n\\ r\leq k_1+k_2}}}\times
P^{k_1}(w_1)P^{k_2}(w_2)x^t.\label{subeq1B}
\end{align}
Since \(P^{k_1}(w_1)P^{k_2}(w_2)=q^{k_1+k_2}p^{k_1m_1+k_2(m+1)}\), from
 (\ref{subeq1A}) and (\ref{subeq1B}), we have
\begin{align*}
&P(\bar{E}_{n,m}(X_1^n)=t)=\\
&\sum_{\substack{k_1,k_2:\\(m+1)k_1+(m+2)k_2\leq n,\\ t\leq k_1+k_2}} 
(-1)^{k_1-t}\dbinom{n-mk_1-(m+1)k_2}{k_1,k_2}\dbinom{k_1+k_2}{t} q^{k_1+k_2}p^{k_1m_1+k_2(m+1)}.
\end{align*}

By similar manner to Lemma~\ref{mainlemma}, we have part 5.
\qed

\section{Algorithm and complexity}\label{seccomplexity}
We study algorithm and complexity to compute (\ref{subexplicitform}).
The basic idea of our algorithm 
is similar to that of bucket sort (\cite{IntroAlgoBook}).
When  \(P(C_{n,(m_1,\ldots,m_h)} >t)\) is negligible for some \(t\), it is suffice to compute \(P(C_n)=s\) for \(s=0,\ldots,t\).

Let \(\Zb_{\geq 0}:=\{0,1,2,\ldots\}\) and 
\begin{align}
G_s:&=\{ (k_1,\ldots,k_h)\in\Zb_{\geq 0}^{h}\mid
 \sum_{1\leq i\leq h}k_i|w_i|\leq n, \sum (i-1)k_i\leq s\leq \sum ik_i\}.\label{defG}
\end{align}
\begin{lemma}
\[\bigcup_{0\leq s\leq t} G_s=\{ (k_1,\ldots,k_h)\in\Zb_{\geq 0}^{h}\mid
 \sum_{1\leq i\leq h}k_i|w_i|\leq n, \sum (i-1)k_i\leq t\}.\]
\end{lemma}
Proof) The relation \(\subseteq\) is clear. Assume that \((k_1,\ldots,k_h)\in \{\sum_{1\leq i\leq h}k_i|w_i|\leq n, \sum (i-1)k_i\leq t\}\).
If \(t\leq \sum ik_i\) then \((k_1,\ldots,k_h)\in G_t\). If \(\sum ik_i<t\) then there is a nonnegative \(s<t\) such that  \((k_1,\ldots,k_h)\in G_s\). We have \((k_1,\ldots,k_h)\in\bigcup_{0\leq s\leq t} G_s\).
\qed

The following Algorithm A compute \(P(C_n)=s\) for all \(s=0,\ldots,t\).

\noindent
{\bf Algorithm A}\\
1.~Initialize \(P(C_n=s)=0\) for all \(s=0,\ldots,t\).\\
2.~Enumerate all nonnegative vectors \((k_1,\ldots,k_h)\in \bigcup_{0\leq s\leq t}G_s\).\\
\phantom{2. }For each vector \((k_1,\ldots,k_h)\in G_s\) and \(s\), set
\begin{align*}
P(C_{n,(w_1,\ldots,w_h)}(X_1^n)=s):=&P(C_{n,(w_1,\ldots,w_h)}(X_1^n)=s)\\
&+(-1)^{\sum ik_i-s}\dbinom{n-\sum |w_i|k_i+\sum k_i }{k_1,\ldots,k_h} \dbinom{\sum k_i}{\sum ik_i-s}\prod P^{k_i}(w_i).
\end{align*}
3.~Output \(P(C_n=s)\) for all \(s=0,\ldots,t\).
\medskip

Since Algorithm A enumerates all combination of \(s,k_1,\ldots,k_h\) such that \((k_1,\ldots,k_h)\in G_s\) for all \(s=0,\ldots,t\) in (\ref{subexplicitform}), 
Algorithm A correctly computes
\(P(C_{n,(w_1,\ldots,w_h)}(X_1^n)=s)\) for all \(s=0,\ldots,t\).

In Theorem~\ref{propcomputationalcomplexity}, we show the size of \(\bigcup_{0\leq s\leq t} G_s\) and the number of operands that are required 
to enumerate \(\bigcup_{0\leq s\leq t} G_s\), which are the bottle necks of time and space complexity of Algorithm A, respectively.
For simplicity, we do not study space and computational complexity as a function of input length.
Algorithm and its time and space complexity for computing \(D_n\) are similar to those of \(C_n\).

Let 
\[F(t,h):=\vert \{ (k_1,\ldots,k_h)\in\Zb_{\geq 0}^h\mid \sum_{1\leq i\leq h}ik_i\leq t\}\vert\text{ and } F(t):=F(t,t),\]
where \(|S|\) is the number of the elements of a finite set \(S\).
\begin{lemma}[Special case of Theorem 4.9.2 in p.~96 \cite{remirez}]\label{lemmaPartitionA}
Fix h. Then
\begin{equation}\label{shormod}
F(t,h)\sim \frac{t^{h}}{(h!)^2}.
\end{equation}
\end{lemma}
Proof) Fix \(h\) in Theorem 4.9.2 in p.~96 \cite{remirez}.\qed
\begin{lemma}[Special case of Meinardus's theorem (Chap.~6 \cite{andrews})]\label{lemmaPartitionB}
\begin{equation}\label{eqnMeinardus}
F(t)\sim \frac{1}{2\pi\sqrt{2t}}\exp(\pi(\frac{2}{3}t)^{\frac{1}{2}})
\end{equation}
\end{lemma}
Proof) Since \(\sum ik_i\leq t\Leftrightarrow\exists s\geq 0\ s+\sum ik_i=t\),
we have \(\sum_n F(n)q^n= \Pi_n (1-q^n)^{-a_n}\) for \(|q|<1\) where \(a_1=2\) and \(a_n=1\) for \(n\geq 2\). 
By applying Theorem 6.2 in Chap.~6 \cite{andrews}, we have the lemma 
(cf.~Theorem 6.3 in \cite{andrews} for the case that \(a_n=1\) for all \(n\)).\qed

For example \(F(30,2)=256\) and  \(F(30)=28629\); corresponding right-hand sides of (\ref{shormod}) and (\ref{eqnMeinardus}) are 225 and
approximately 25967, respectively.
\begin{theorem}\label{propcomputationalcomplexity}
For given \(n,h,t,\) and \(w_1,\ldots,w_h\), \\
1. Fix \(t\). Then \(|\bigcup_{1\leq s\leq t}G_s|=O(n)\).\\
2. Fix \(h\). Then \(|\bigcup_{1\leq s\leq t}G_s|=O(nt^{h-1})\).\\
3. Let \(\alpha>0\) and  \(t:=\frac{3\alpha^2}{2\pi^2}(\log n)^2\).
Then \(|\bigcup_{1\leq s\leq t}G_s|=O(\frac{1}{2\sqrt{3}\alpha |w_1|\log n}n^{1+\alpha})\).\\
4. The number of operands required to enumerate \(\bigcup_{0\leq s\leq t} G_s\) is \(O(t)\).
\end{theorem}
Proof)
Since 
\(
\bigcup_{0\leq s\leq t} G_s\subseteq \{ (k_1,\ldots,k_h)\mid k_1\leq\lfloor \frac{n}{|w_1|}\rfloor, \sum (i-1)k_i\leq t\},
\)
 we have
\(
|\bigcup_{0\leq s\leq t} G_s|\leq \frac{n}{|w_1|} F(t,h-1).
\)
From Lemma~\ref{lemmaPartitionA} and \ref{lemmaPartitionB}, we have Part 1--3.
Part 4 is clear. 
\qed
\begin{remark}
To compute exact distributions by Markov imbedding method \cite{Fu94}, we need to calculate \(M^n\) for  sample size \(n\) and
\(m\times m\) matrix \(M\) with \(m=O(n)\).
The  number of arithmetical operations to compute \(M^2\) is \(O(n^{2.81})\) with Strassen algorithm (\cite{IntroAlgoBook}) and those of \(M^n\) is 
\(O(n^{2.81}\log n)\). The number of operands to compute Markov imbedding method with Strassen algorithm is \(O(n^2)\), see \cite{Bailey90}.
\end{remark}
Let
\[G(n,s_1,\ldots, s_h):=|\{(k_1,\ldots,k_h)\vert s_1\leq k_1,\ldots,s_h\leq k_h, \sum |w_i|k_i\leq n\}|.\]
The quantity \(G\) is the bottle neck of time complexity to compute (\ref{eqmain}).
\begin{theorem}
\[
G(n,s_1,\ldots, s_h)=
\begin{cases}
F(n-\sum s_i|w_i|, h) &\text{if } \sum |w_i|s_i\leq n,\\
0 &\text{else.}
\end{cases}
\]
If \(|w_i|=m\) for \(i=0,\ldots, h\) then
\[G(n,s_1,\ldots, s_h)=\binom{ \lfloor\frac{n}{m}\rfloor -\sum s_i +h-1}{h-1}\text{ if }\sum s_i\leq \lfloor\frac{n}{m}\rfloor.\]
\end{theorem}
Proof)
Assume that \(s_1\leq k_1,\ldots,s_h\leq k_h\).
Then \(\sum |w_i|k_i\leq n\Leftrightarrow \sum |w_i| (k_i-s_i)\leq n-\sum |w_i|s_i\) and we have the first equation.
Since \(\sum |w_i| (k_i-s_i)\leq n-\sum |w_i|s_i\Leftrightarrow \sum (k_i-s_i) \leq \lfloor \frac{n}{m}\rfloor-\sum s_i\) 
if \(|w_i|=m\) for all \(i\), by occupancy theorem (5.2) in p.38 \cite{feller1}, we have the lemma.\qed

\section{Distance of distributions}\label{Secdist}
We show that \(P(C_{n,(w_1,\ldots,w_d)}=t)\) and \(P(D_{n,(m_1,\ldots,m_h)}=t)\) uniformly
converge to  \(P(C_{n,(w_1,\ldots,w_h)}=t)\) and \(P(D_{n,(m_1,\ldots,m_h)}=t)\) as \(d\to h\), respectively.
\begin{prop}\label{diffDistA} Let \(X_1,\ldots,X_n\) be i.i.d.~binary random variables from \(P(X_i=0)=p\). Assume \(d<h\). Then 
\begin{align*}
&\sup_t | P( C_{n,(w_1,\ldots,w_d)}(X_1^n)=t) -P( C_{n,(w_1,\ldots,w_h)}(X_1^n)=t) | \leq (n-|w_{d+1}|+1)P(w_{d+1})\text{ and }\\
&\sup_t | P(D_{n,(m_1,\ldots,m_d)}(X_1^n)=t) -P(D_{n,(m_1,\ldots,m_h)}(X_1^n)=t)|\leq 2(n+1-m_{d+1})p^{m_{d+1}}.
\end{align*}
\end{prop}
Proof)
Assume that \(\Nmb^{\prime}(w_1,\ldots,w_h; x_1^n)=(k_1,k_2,\ldots,k_h)\).
By (\ref{defCnfun}), \(C_{n,(w_1,\ldots,w_d)}(x_1^n)=C_{n,(w_1,\ldots,w_h)}(x_1^n)\) if \(k_{d+1}=\cdots=k_h=0\).
Then for all \(t\),
\begin{align}
&| P( C_{n,(w_1,\ldots,w_d)}(X_1^n)=t) -P( C_{n,(w_1,\ldots,w_h)}(X_1^n)=t) |\nonumber\\
&\leq P\{\text{ there is }h\text{ such that }d+1\leq h\text{ and }k_h>0\}\nonumber\\
&\leq (n-|w_{d+1}|+1)P(w_{d+1}).\label{approxUpperA}
\end{align}
Let \(w_i=10^{m_i}\) for all \(i\).
By Theorem~\ref{mainth2}, for all \(t\),
\begin{align*}
& | P(D_{n,(m_1,\ldots,m_d)}(X_1^n)=t) -P(D_{n,(m_1,\ldots,m_h)}(X_1^n)=t)|\\
& \leq \frac{1}{q} \vert P(C_{n+1,(w_1,\ldots,w_d)}=t)-P(C_{n+1,(w_1,\ldots,w_h)}=t)\vert +\frac{p}{q}
\vert P(C_{n,(w_1,\ldots,w_d)}=t)-P(C_{n,(w_1,\ldots,w_h)}=t)\vert \\
&\leq 2(n+1-m_{d+1})p^{m_{d+1}},
\end{align*}
where the last inequality follows from (\ref{approxUpperA}) and \(P(w_{d+1})=qp^{m_{d+1}}\).
\qed

\begin{figure}[h]
\begin{center}
\includegraphics[width=10cm]{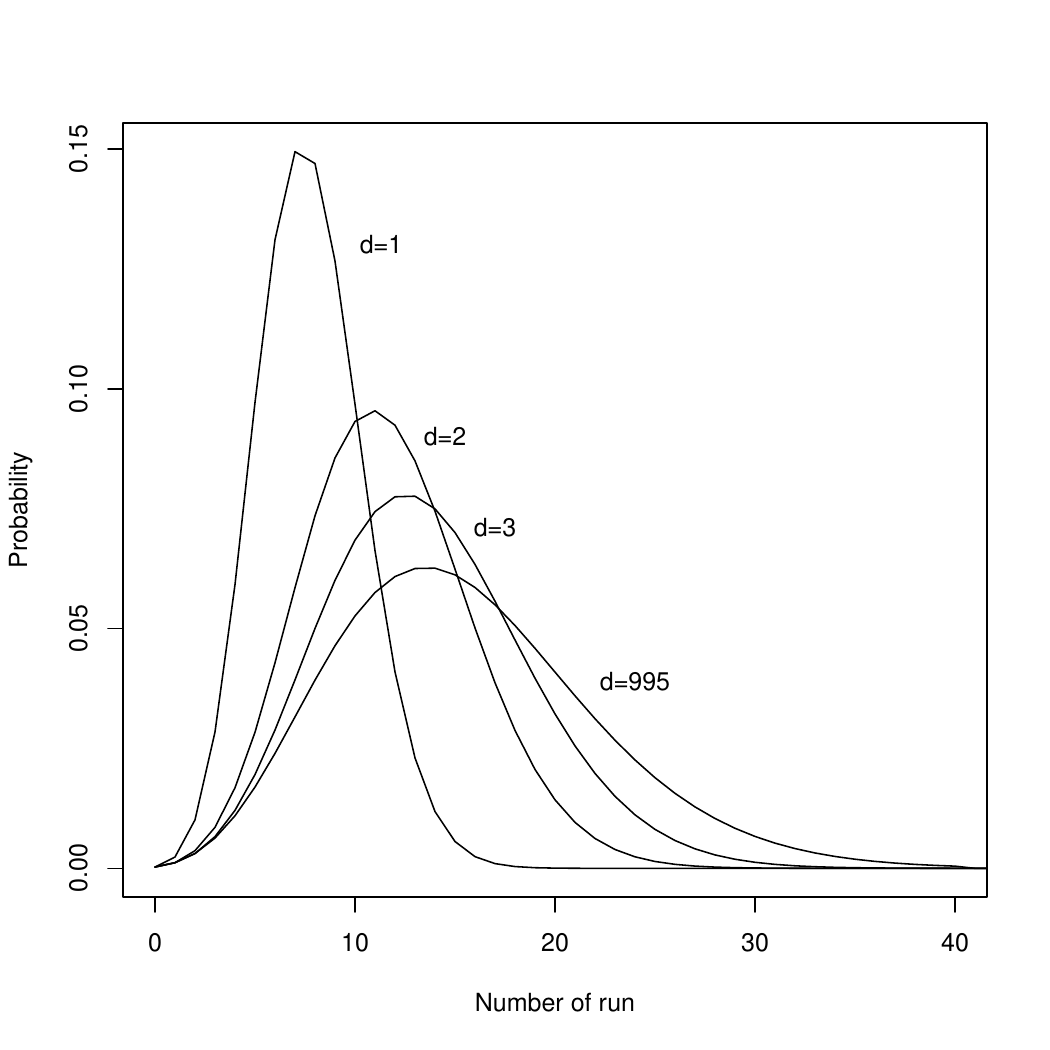}
\caption{Graph of distributions
}\label{figA}
\end{center}
\end{figure}

Assume that  \(X_1,\ldots,X_n\) be i.i.d.~binary random variables from \(P(X_i=0)=0.5\).
Let 
\[
\dist(d,h\mid r):=\sup_{0\leq t\leq r} | P(D_{n,(m_1,\ldots,m_d)}(X_1^n)=t) -P(D_{n,(m_1,\ldots,m_h)}(X_1^n)=t)|.
\]
Table~\ref{tabA} shows  numerical calculations of \(\dist(d,995|40)\) for \(n=1000, d=1,3,5,7,9\), and \(m_i=5+i\) for \(i=1,2,\ldots,995\).
Figure~\ref{figA} shows graphs of \(P(D_{n,(m_1,\ldots,m_d)}(X_1^n)=t)\) for \(d=1,2,3\), and \(995\).
\begin{table}[h]
\begin{center}
\caption{
Distance of distributions
}\label{tabA}
\begin{tabular}{lccccc}
\(d\)  & 1 & 3 & 5 & 7 & 9 \\ \hline
\(\dist(d,995|40)\)  & 0.117859 & 0.0168652 & 0.0036909 & 0.0009005 & 0.0002248 
\end{tabular}
\end{center}
\end{table}

\begin{center}
{\bf Acknowledgement}
\end{center}
This work was supported by the Research Institute for Mathematical Sciences,
an International Joint Usage/Research Center located in Kyoto University.
The author thanks Prof.~Shigeki Akiyama (Tsukuba Univ.) for discussions.


{\small
\begin{thebibliography}{10}

\bibitem{akihirano2000}
S.~Aki and K.~Hirano.
\newblock Numbers of success-runs of specified length until certain stopping
  time rules and generalized binomial distributions of order k.
\newblock {\em Ann.~Inst.~Statist.~Math.}, 52(4):767--777, 2000.

\bibitem{hirano84}
S.~Aki, H.~Kuboki, and K.~Hirano.
\newblock On discrete distributions of order k.
\newblock {\em Ann.~Inst.~Statist.~Math.}, 36:431--440, 1984.

\bibitem{remirez}
J.~L.~Ram\'{\i}rez Alfons\'{\i}n.
\newblock {\em The {D}iophantine {F}robenius {P}roblem}.
\newblock Oxford University Press, 2005.

\bibitem{andrews}
G.~E. Andrews.
\newblock {\em The theory of partitions}.
\newblock Cambridge University Press, 1984.

\bibitem{antzou99}
D.~L. Antzoulakos and S.~Chadjiconstantindis.
\newblock Distributions of numbers of success runs of fixed length in {M}arkov
  dependent trials.
\newblock {\em Ann.~Inst.~Statist.~Math.}, 53(3):599--619, 2001.

\bibitem{Bailey90}
D.~H. Bailey, K.~Lee, and H.~D. Simon.
\newblock Using {S}trassen's algorithm to accelerate the solution of linear
  systems.
\newblock {\em Journal of Supercomputing}, 4:357--371, 1990.

\bibitem{koutrasBook}
N.~Balakrishnan and M.~V. Koutras.
\newblock {\em Runs and scans with applications}.
\newblock John Wiley \& Sons, 2002.

\bibitem{Bassinoetal2011}
F.~Bassino, J.~Cl{\'e}ment, and P.~Micod{\`e}me.
\newblock Counting occurrences for a finite set of words: combinatorial
  methods.
\newblock {\em ACM Trans.~Algorithms.}, 9(4):Article No.~31, 2010.

\bibitem{bertheRigo}
V.~Berth\'{e} and M.~Rigo.
\newblock {\em Combinatorics, words and symbolic dynamics}.
\newblock Encyclopedia of Mathematics and Its Applications 159. Cambridge
  University Press, 2016.

\bibitem{blomThorburn}
G.~Blom and D.~Thorburn.
\newblock How many random digits are required until given sequences are
  obtained?
\newblock {\em J.~Appl.~Probab.}, 19(3):518--531, 1982.

\bibitem{chrysaphinou88}
O.~Chrysaphinou and S.~Papastavridis.
\newblock A limit theorem on the number of overlapping appearances of a pattern
  in a sequence of independent trials.
\newblock {\em Probab.~Theory Related Fields}, 79:129--143, 1988.

\bibitem{IntroAlgoBook}
T.~H. Cormen, C.~E. Leiseron, R.~L. Rivest, and C.~Stein.
\newblock {\em Introduction to algorithms}.
\newblock MIT Press, 3rd edition, 2009.

\bibitem{feller1}
W.~Feller.
\newblock {\em An Introduction to probability theory and its applications
  {V}ol.~1}.
\newblock Wiley, 3rd revised edition, 1970.

\bibitem{Flajolet2009}
P.~Flajolet and R.~Sedgewick.
\newblock {\em Analytic Combinatorics}.
\newblock Cambridge University Press, 2009.

\bibitem{Fu94}
J.~C. Fu and M.~V. Koutras.
\newblock Distribution theory of runs: a {M}arkov chain approach.
\newblock {\em J.~Amer.~Statist.~Assoc.}, 89(427):1050--1058, 1994.

\bibitem{Godbole90}
A.~P. Godbole.
\newblock Specific formulae for some success run distributions.
\newblock {\em Statist.~Probab.~Lett.}, 10:119--124, 1990.

\bibitem{Godbole92}
A.~P. Godbole.
\newblock The exact and asymptotic distribution of overlapping success runs.
\newblock {\em Comm.~Statist.~Theory Methods}, 21:953--967, 1992.

\bibitem{Goulden83}
I.~Goulden and D.~Jackson.
\newblock {\em Combinatorial Enumeration}.
\newblock John Wiley, 1983.

\bibitem{Guibas81}
L.~Guibas and A.~Odlyzko.
\newblock String overlaps, pattern matching, and nontransitive games.
\newblock {\em J.~Combin.~Theory Ser.~A}, 30:183--208, 1981.

\bibitem{hanAki}
S.~Han and S.~Aki.
\newblock A unified approach to binomial-type distributions of order k.
\newblock {\em Commun.~Statist.~Theor.~Meth.}, 29:1929--1943, 2000.

\bibitem{Hirano86}
K.~Hirano.
\newblock Some properties of the distributions of order k.
\newblock pages 43--53, 1986.
\newblock Fibonacci Numbers and their Applications, A.~N.~Phillipou,
  A.~F.~Horadam and G.~E.~Bergum eds, Reidel.

\bibitem{JaquetandSzpankowski}
P.~Jacquet and W.~Szpankowski.
\newblock {\em Analytic Pattern Matching}.
\newblock Cambridge University Press, 2015.

\bibitem{koutras97}
M.~V. Koutras and V.~A. Alexandrou.
\newblock Non-parametric randomness tests based on success runs of fixed
  length.
\newblock {\em Statist.~Probab.~Lett.}, 32:393--404, 1997.

\bibitem{ling88}
K.~D. Ling.
\newblock On binomial distributions of order k.
\newblock {\em Statist.~Probab.~Letters}, 6:247--250, 1988.

\bibitem{lothire}
M.~Lothaire.
\newblock {\em Applied Combinatorics on words}.
\newblock Encyclopedia of Mathematics and Its Applications 105. Cambridge
  University Press, 2005.

\bibitem{Makri2007}
F.~S. Makri, A.~N. Philippou, and Z.~M. Psillakis.
\newblock Shortest and longest length of success runs in binary sequences.
\newblock {\em J.~Statist.~Plan.~Inference}, 137:2226--2239, 2007.

\bibitem{Makri2015}
F.~S. Makri and Z.~M. Psillakis.
\newblock On l-overlapping runs of ones of length k in sequences of independent
  binary random vairables.
\newblock {\em Commun.~Statist.~Theor.~Meth.}, 44:3865--3884, 2015.

\bibitem{Mood40}
A.~M. Mood.
\newblock The distribution theory of runs.
\newblock {\em Ann.~Math.~Statist}, 11(4):367--392, 1940.

\bibitem{Muselli96}
M.~Muselli.
\newblock Simple expressions for success run distributions in {B}ernoulli
  trials.
\newblock {\em Statist.~Probab.~Lett.}, 31:121--128, 1996.

\bibitem{phili83}
A.~N. Philippou, C.~Georghiou, and G.~N. Philippou.
\newblock A generalized geometric distribution and some of its properties.
\newblock {\em Statist.~Probab.~Letters}, 1:171--175, 1983.

\bibitem{Phillipou86}
A.~N. Phillipou and F.~S. Makri.
\newblock Successes, runs and longest runs.
\newblock {\em Statist.~Probab.~Lett.}, 4:211--215, 1986.

\bibitem{RegnierSpankowski98}
M.~R{\'e}gnier and W.~Szpankowski.
\newblock On pattern frequency occurrences in a {M}arkovian sequence.
\newblock {\em Algorithmica}, 22(4):631--649, 1998.

\bibitem{Riordan}
J.~Riordan.
\newblock {\em Introduction to combinatorial analysis}.
\newblock John Wiley, 1958.

\bibitem{bookRobin}
S.~Robin, F.~Rodolphe, and S.~Schbath.
\newblock {\em DNA, words and models}.
\newblock Cambridge University Press, english edition, 2005.

\bibitem{rukhin08}
A.~L. Rukhin and Z.~Volkovich.
\newblock Testing randomness via aperiodic words.
\newblock {\em J.~Stat.~Comput.~Simul.}, 78(12):1133--1144, 2008.

\bibitem{takahashi2022ProbSympo}
H.~Takahashi.
\newblock Some explicit formulae for the distributions of words.
\newblock {\em RIMS K{\^o}ky{\^u}roku, Kyoto University}, 2246:16--23, Apl.
  2023.
\newblock Probability Symposium 2022.

\bibitem{uppuuri}
V.~R.~R. Uppuluri and S.~A. Patil.
\newblock Waiting times and generalized {F}ibonacci sequences.
\newblock {\em Fibonacci Quart.}, 21:242--249, 1983.

\bibitem{Wald40}
A.~Wald and J.~Wolfowitz.
\newblock On a test whether two samples are from the same population.
\newblock {\em Ann.~Math.~Statist}, 11(2):147--162, 1940.

\bibitem{waterman}
M.~S. Waterman.
\newblock {\em Introduction to computational biology}.
\newblock Chapman \& Hall, New York, 1995.

\bibitem{wolf88}
E.~Z. Zehavi and J.~K. Wolf.
\newblock On runlength codes.
\newblock {\em IEEE Trans.~Inform.~Theory}, 34(1):45--53, 1988.

\end{thebibliography}
}
\end{document}